%% file: main.tex
\documentclass[
    aps,
    pra,
    reprint,
    superscriptaddress,
]{revtex4-2}
\usepackage[dvipdfmx]{graphicx}
\usepackage{hyperref}
\usepackage{amsmath}
\usepackage{amsfonts}
\usepackage{amssymb}
\usepackage{mathtools}
\usepackage{amsthm}
\usepackage{color}
\usepackage{enumitem}

\begin{document}
\title{Programming sequential deployment of origami via kinematic transition fronts}
\author{Rinki Imada}
\email{imada.rinki@jaxa.jp}
\affiliation{
    Institute of Space and Astronautical Science,
    Japan Aerospace Exploration Agency, 
    Sagamihara-shi,
    Kanagawa 252-5210,
    Japan
}
\affiliation{
    Department of General Systems Studies, 
    Graduate School of Arts and Sciences, 
    The University of Tokyo,
    Meguro-ku,
    Tokyo 153-8902,
    Japan
}
\author{Tomohiro Tachi}
\email{tachi@idea.c.u-tokyo.ac.jp}
\affiliation{
    Department of General Systems Studies, 
    Graduate School of Arts and Sciences, 
    The University of Tokyo,
    Meguro-ku,
    Tokyo 153-8902,
    Japan
}

\begin{abstract}
    Propagating transition fronts, in which local interactions sequentially trigger state changes, are widely observed across natural, biological, and engineered systems.
    While such propagation has been engineered using energy-driven instabilities, front propagation governed purely by geometric constraints remains underexplored and lacks a general design framework. 
    In particular, how to program sequential deployment in origami through such kinematic propagation remains an open challenge.
    Here, we develop a systematic design framework for kinematic transition fronts based on their correspondence with heteroclinic orbits in discrete dynamical systems.
    Focusing on strips of developable and flat-foldable degree-4 origami vertices, we show that asymmetric coupling between adjacent creases produces nonlinear recurrence relations whose composition generically gives rise to heteroclinic orbits connecting developed and flat-folded states, enabling domino-like sequential deployment.
    We further show that macroscopic shape can be programmed independently of propagation behavior by exploiting invariances in the recurrence relation, and illustrate the approach through a representative thick-panel origami prototype.
    These results enable programmable sequential deployment in origami via transition fronts, while also establishing a general framework for kinematic transition fronts in geometrically constrained systems.
\end{abstract}

\maketitle

\section{Introduction}
\begin{figure*}[t]
    \centering
    \includegraphics[width=\linewidth]{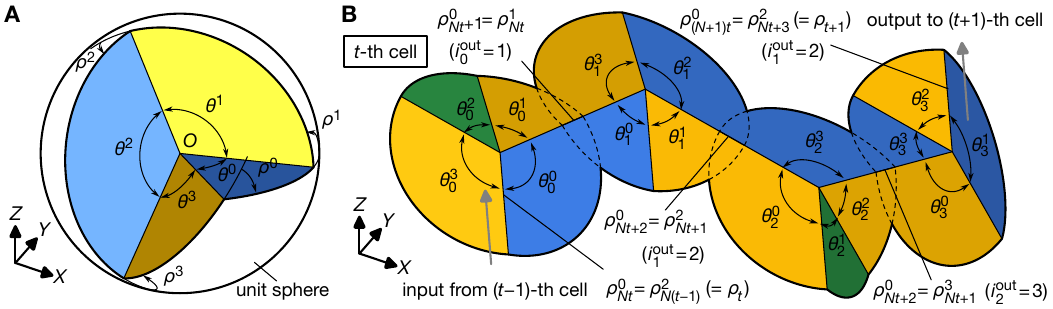}
    \caption{
        (A) Parameterization of the design and kinematics of a single degree-4 vertex by the sector angles $(\theta^{i})_{i=0,1,2,3}$ and fold angles $(\rho^{i})_{i=0,1,2,3}$.
        The sector angle $\theta^i$ denotes the angle between the $i$-th and $(i+1)$-th creases, and the fold angle $\rho^i$ denotes the exterior dihedral angle between the two facets sharing the $i$-th crease.
        A valley (mountain) fold corresponds to $\rho^i>0$ ($\rho^i<0$).
        A single degree-4 vertex can be represented as a spherical quadrilateral with fixed arc lengths $(\theta^i)_{i=0,1,2,3}$ by intersecting the crease pattern with a unit sphere centered at the vertex $O$, i.e., a spherical 4R mechanism.
        (B) Example of the unit-cell of a periodic degree-4 origami strip ($N=4, (i^{\mathrm{out}}_{0},i^{\mathrm{out}}_{1},i^{\mathrm{out}}_{2},i^{\mathrm{out}}_{3})=(1,2,3,2)$).
    }
    \label{fig:parameterization}
\end{figure*}
Domino-like transitions (i.e., propagating transition fronts), in which local interactions trigger state changes, are a ubiquitous propagation phenomenon across natural, biological, and engineered systems.
Examples include the cascading release of mechanical instabilities~\cite{nadkarni2014dynamics,raney2016stable}, calcium waves in biological tissues~\cite{leybaert2012intercellular}, perversion points in helical structures~\cite{goriely1998tendril,mcmillen2002tendril}, and localized flipping in a ribbon-spread deck of cards~\cite{gardner2002playing}.
In recent years, mechanical metamaterials have provided a platform for engineering such propagating transition fronts by designing interactions between constituent elements~\cite{kochmann2017exploiting,deng2021nonlinear,veenstra2024non,janbaz2024diffusive}.
A representative example is the transition wave, where the release of stored elastic energy in bistable units triggers sequential switching in neighboring elements~\cite{nadkarni2016unidirectional,jin2020guided,zareei2020harnessing,yasuda2020transition}.
Propagating transition fronts can also arise from geometric interactions, without relying on energy barriers, as a purely kinematic, nonlinear phenomenon.
Such phenomena have been reported in systems such as the Kane-Lubensky chain~\cite{kane2014topological} in the context of topological mechanics~\cite{chen2014nonlinear,lubensky2015phonons,chen2016topological,rocklin2017transformable,mao2018maxwell}, where geometric constraints dictate the order and direction of transitions.
These mechanisms are not limited to idealized linkage models but can be realized in a broader class of systems governed by geometric compatibility.
However, such kinematic transition fronts remain largely underexplored and lack a general design framework.

Origami provides a particularly compelling platform in this regard, as it offers a geometric principle for compactly storing large systems.
While many canonical mechanisms, such as the Miura-ori, require globally coordinated actuation in which all creases must be driven simultaneously, connecting a compact stowed state and a deployed state via a propagating transition front would enable deployment through spatially localized, sequential actuation.
This removes the need for global synchronization and facilitates scalable actuation strategies.
Moreover, such sequential deployment can reduce the swept volume during actuation, since only a localized region undergoes large motion at each stage, which is critical in spatially confined environments.
These features are advantageous for applications such as minimally invasive medical devices (e.g., catheters and endoscopes), inspection tools for confined infrastructures, and the installation of structures within existing built environments.
It may also simplify the programming of folding sequences by embedding actuation order into the geometry, avoiding the need for crease-by-crease stimulus tuning in conventional self-folding approaches~\cite{liu2017sequential}.
However, despite these potential advantages, a general design framework for sequential deployment in origami remains lacking.

Here, we present a systematic design strategy for origami-based mechanisms that realize kinematic propagation of transition fronts.
Specifically, we focus on a one-degree-of-freedom strip of connected degree-4 vertices, a building block in origami kinematics.
Although origami structures composed of degree-4 vertices, such as the Miura-ori, have been extensively studied in the contexts of kinematics~\cite{tachi2009generalization,tachi2011one,lang2017twists,hull2020origametry}, discrete differential geometry~\cite{schief2008integrability,stachel2010kinematic,izmestiev2017classification,sharifmoghaddam2023generalizing,he2026infinitely}, and engineering~\cite{schenk2013geometry,filipov2015origami,pratapa2019geometric,dudte2021additive}, the realization of propagating transition fronts in such systems has remained elusive and was previously considered unattainable~\cite{chen2016topological}. 
Recent studies, however, have suggested that this limitation can be overcome by modifying the connectivity of crease patterns~\cite{imada2025kinematic}.
The degree-4 origami strip is a one-dimensional Maxwell lattice~\cite{lubensky2015phonons,mao2018maxwell,imada2025maxwell}, in which specifying the folded state of the boundary vertex determines the states of subsequent vertices in the bulk~\cite{chen2016topological,imada2025kinematic}.
Previous work has shown that, when the recurrence relation between adjacent cells is interpreted as a discrete dynamical system~\cite{imada2022geometry,imada2023undulations,imada2025maxwell,imada2025kinematic,imada2026nonlinear,kim2022nonlinear,zhou2017kink}, heteroclinic orbits emerging in sufficiently nonlinear regimes give rise to propagating transition fronts~\cite{imada2025kinematic}.
However, only a few specific examples have been reported, and a general design methodology has not yet been established; moreover, these findings have remained largely theoretical without experimental validation.

Building on the correspondence between heteroclinic orbits in discrete dynamical systems and the kinematic propagation of transition fronts, we design degree-4 origami strips that deploy in a domino-like manner and validate them using physical prototypes.
We focus on strips composed of periodically connected, developable and flat-foldable degree-4 vertices.
Based on the asymmetry in the folding angles of adjacent creases at each vertex, we show that, when adjacent creases are used to connect vertices, heteroclinic orbits generically emerge over a broad parameter range, leading to domino-like deployment, as demonstrated through representative design examples.
Furthermore, we show that, by transforming the crease pattern under conditions that preserve the recurrence relation, the macroscopic shape of the deployed and flat-folded configurations can be programmed independently while maintaining the propagation behavior.
In addition, we implement a representative design as a thick-panel origami structure using a thickness-accommodation technique, and experimentally observe the propagation of the transition in a 3D-printed prototype.
These results not only enable the systematic design of sequential deployment in origami, but also provide a general framework for engineering kinematic transition fronts beyond origami, with potential applications in deployable structures, robotics, mechanical information propagation, and metamaterials.

\section*{Results}
\begin{figure*}[t]
    \centering
    \includegraphics[width=\linewidth]{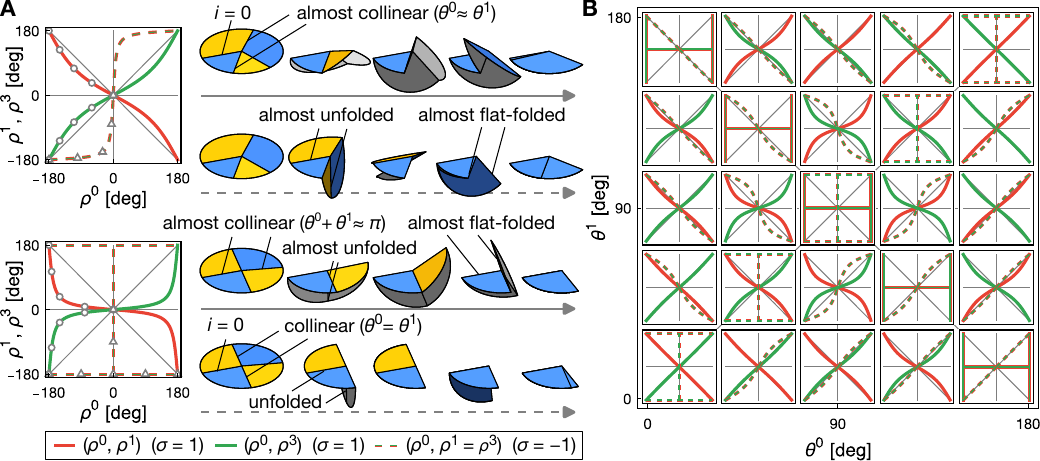}
    \caption{
        (A) Visualization of the relationship between adjacent fold angles and the folding motions for two representative examples (top: $(\theta^{0},\theta^{1})=(60^\circ,55^\circ)$; bottom: $(\theta^{0},\theta^{1})=(85^\circ,85^\circ)$).
        In the left panels, the folding modes corresponding to $\sigma=1$ ($\sigma=-1$) are shown as solid (dashed) curves in the parameter space.
        In the right panels, folding motions from the developed state to the flat-folded state are visualized by sampling configurations along these curves, with circles (triangles) indicating $\sigma=1$ ($\sigma=-1$).
        Both examples have nearly collinear creases, causing the highly asymmetric folding behavior.
        In the bottom example, a pair of creases ($i=1,3$) is exactly collinear ($\theta^{0}=\theta^{1}$), resulting in a decoupling between the two pairs $(i=0,2)$ and $(i=1,3)$ when $\sigma=-1$.
        In this case, the collinear pair can be flat-folded first without folding the other pair, corresponding to the vertical branch; upon reaching the flat-folded state $((\rho^{0},\rho^{1})=(0,-\pi)$ or $(0,\pi))$, the other pair $(i=0,2)$ subsequently aligns and can also be flat-folded, corresponding to the horizontal branches in the left panel.
        (B) Dependence of the adjacent fold-angle relations on $(\theta^{0},\theta^{1})$.
        For representative sector angles, the relations are visualized as curves in the same format as in (A) and arranged over the $(\theta^{0},\theta^{1})$-space.
        For panels along the diagonals $\theta^{0}=\theta^{1}$ ($\theta^{0}+\theta^{1}=\pi$), the curves for $\sigma=-1$ ($\sigma=1$) show the decoupling of the two pairs of opposite creases, similar to the behavior shown in the bottom panels of (A).
        Except for the singular cases, each curve intersects the diagonal lines $\rho^{1}=\pm\rho^{0}$ and $\rho^{3}=\pm\rho^{0}$ only at $(0,0)$, $(\pm\pi,\pm\pi)$, or $(\pm\pi,\mp\pi)$ and the slope at these intersection points depends on $(\theta^{0},\theta^{1})$.
        This reflects that the cosines $\cos\rho^{0}$ and $\cos\rho^{1}=\cos\rho^{3}$ are related by a symmetric fractional linear transformation, for which $\cos\rho^{0}=\pm1$ are the only fixed points.
    }
    \label{fig:asymmetry}
\end{figure*}
\subsection*{Parameterization}
We introduce a minimal parameterization of the crease pattern design and kinematics of periodic origami strips, following our previous work~\cite{imada2025kinematic}.
A single degree-4 vertex is parameterized by labeling its creases counterclockwise as $i=0,1,2,3$, and introducing the sector angles $(\theta^{0}, \theta^{1}, \theta^{2}, \theta^{3})$ and the fold angles $(\rho^{0}, \rho^{1}, \rho^{2}, \rho^{3})$, where $\theta^{i}\in(0,\pi)$ and $\rho^{i}\in [-\pi,\pi]$ (Fig.~\ref{fig:parameterization}(A)).
Imposing developability and flat-foldability leads to the constraints $\theta^{0}+\theta^{1}+\theta^{2}+\theta^{3}=2\pi$ and $\theta^{0}-\theta^{1}+\theta^{2}-\theta^{3}=0$, which are equivalent to $\theta^{2}=\pi-\theta^{0}$ and $\theta^{3}=\pi-\theta^{1}$.
For a strip, we label the vertices by $n=0,1,\dots$ and denote the sector angles and fold angles at the $n$-th vertex by $(\theta_{n}^{i})_{i=0,1,2,3}$ and $(\rho_{n}^{i})_{i=0,1,2,3}$.
When vertices are connected in order of $n$, each vertex shares creases with its neighbors, defining input and output creases.
Here, we fix the input crease to $i=0$ and denote the output crease index at the $n$-th vertex by $i_n^{\mathrm{out}} \in \{1,2,3\}$.
The crease pattern of the strip is thus specified by $((\theta_{n}^{i})_{i=0,1,2,3})_{n=0,1,\dots}$ and $(i_n^{\mathrm{out}})_{n=0,1,\dots}$; in this study, we focus on periodic strips with period $N\in\mathbb{Z}_{\geq 1}$, i.e., $\theta_{n+N}^{i}=\theta_{n}^{i}$ and $i_{n+N}^{\mathrm{out}}=i_{n}^{\mathrm{out}}$.
Owing to the periodicity, we define the $t$-th cell by $Nt$-th to $(N(t+1)-1)$-th vertices ($t=0,1,\dots$).
Consequently, a periodic, developable, and flat-foldable strip is fully characterized by the period $N$, the sector angles and connectivity of a unit cell, $((\theta_{n}^{i})_{i=0,1})_{n=0,1,\dots,N-1}$ and $(i_n^{\mathrm{out}})_{n=0,\dots,N-1}$ (Fig.~\ref{fig:parameterization}(B)).
Although crease lengths are additional design parameters, they do not affect the resulting discrete dynamical system and are therefore omitted, except when specifying concrete examples later.

\subsection*{Asymmetric coupling of adjacent fold angles}
The sequence of fold angles in the strip follows the kinematic relations for a developable and flat-foldable degree-4 vertex:
\begin{align}\label{eq:fold-angle}
    \begin{split}
        \rho^{1}=&-\mathrm{sgn}(\rho^{0})\mathrm{sgn}(\cos{\theta^{0}}+\sigma\cos{\theta^{1}})\\
        \quad&\times\arccos\left(\dfrac{A\cos\rho^{0}+B}{B\cos\rho^{0}+A}\right),\\
        \rho^{2}=& \sigma\rho^{0},\\
        \rho^{3}=& -\sigma\rho^{1},\\
        \text{where}\qquad A\coloneq& \cos{\theta^{0}}\cos{\theta^{1}}+\sigma,\, B\coloneq \sin{\theta^{0}}\sin{\theta^{1}},\\
        \sigma \coloneq& \mathrm{sgn}(\rho^{0})\mathrm{sgn}(\rho^{2}).
    \end{split}
\end{align}
The sign parameter $\sigma\in \{-1,1\}$ specifies whether the MV-assignments of opposite fold angles $\rho^0$ and $\rho^2$ are the same or different, which characterizes the folding mode.
Note that Eq.~\eqref{eq:fold-angle} does not apply in the singular cases $(\theta^{0}=\theta^{1},\, \sigma=-1)$ and $(\theta^{0}+\theta^{1}=\pi,\, \sigma=1)$.
Eq.~\eqref{eq:fold-angle} can be derived from spherical trigonometry and is often expressed using alternative variables such as $\tan(\rho^i/2)$~\cite{hull2020origametry,foschi2022explicit}.
From Eq.~\eqref{eq:fold-angle}, the input crease ($i=0$) and the opposite crease ($i=2$) evolve at the same rate; in contrast, the relations between the input crease and the adjacent creases ($i=1,3$) are nonlinear.
Fig.~\ref{fig:asymmetry}(A) and (B) illustrate the adjacent fold-angle relations and the corresponding folding motions for two representative values of $(\theta^{0},\theta^{1})$, as well as their dependence on $(\theta^{0},\theta^{1})$.
As shown in Fig.~\ref{fig:asymmetry}(A), the adjacent creases exhibit asymmetric folding rates: for $\sigma=1$ ($\sigma=-1$), they evolve more slowly (more rapidly) than the input crease when $\rho^{0}$ increases from $0$ to $\pi$ (or $-\pi$), catching up near the flat-folded state.
This asymmetry is quantified by the slopes $|d\rho^{1}/d\rho^{0}|=|d\rho^{3}/d\rho^{0}|$ evaluated at $\rho^{0}=0$ or $\rho^{0}=\pm\pi$, which are reciprocal to each other and given by $\sqrt{(\sigma+\cos(\theta^{0}+\theta^{1}))/(\sigma+\cos(\theta^{0}-\theta^{1}))}=\sqrt{(A-B)/(A+B)}$ and its inverse, respectively.
These quantities are known as the folding multiplier~\cite{schief2008integrability,evans2015twist,evans2015gadgets,tachi2017self}, which has primarily been used to assess rigid foldability, but has not been interpreted as a measure of asymmetry.
As these ratios approach $0$ or $\infty$, the asymmetry becomes increasingly pronounced.
Geometrically, the singular cases ($\theta^0=\theta^1$, or $\theta^0+\theta^1=\pi$) correspond to crease patterns in which opposite creases ($i=1$ and $3$, or $i=0$ and $2$) become collinear in the developed state.
As shown in the bottom panel of Fig.~\ref{fig:asymmetry}(A), when the opposite creases are collinear, they can fold flat without actuating the other creases, i.e., the two pairs of opposite creases become kinematically decoupled.
The asymmetry becomes most pronounced when the opposite creases are nearly collinear so that the two pairs remain only marginally kinematically coupled.
The sequence of input fold angles in the strip, $(\rho_{n}^{0})_{n=0,1,\dots}$, is generated by successive applications of local maps $f_{n}$ to $\rho_{0}^{0}$, where each $f_{n}$ maps the $n$-th input fold angle $\rho_{n}^{0}$ to the output fold angle specified by the index $i_{n}^{\mathrm{out}}$.
By connecting vertices through adjacent creases, the asymmetric response of adjacent fold angles accumulates along the strip, providing the mechanism for domino-like deployment.

\subsection*{Heteroclinic orbits as a design principle for domino-like deployment}
\begin{figure*}[p]
    \centering
    \includegraphics[keepaspectratio,width=\linewidth]{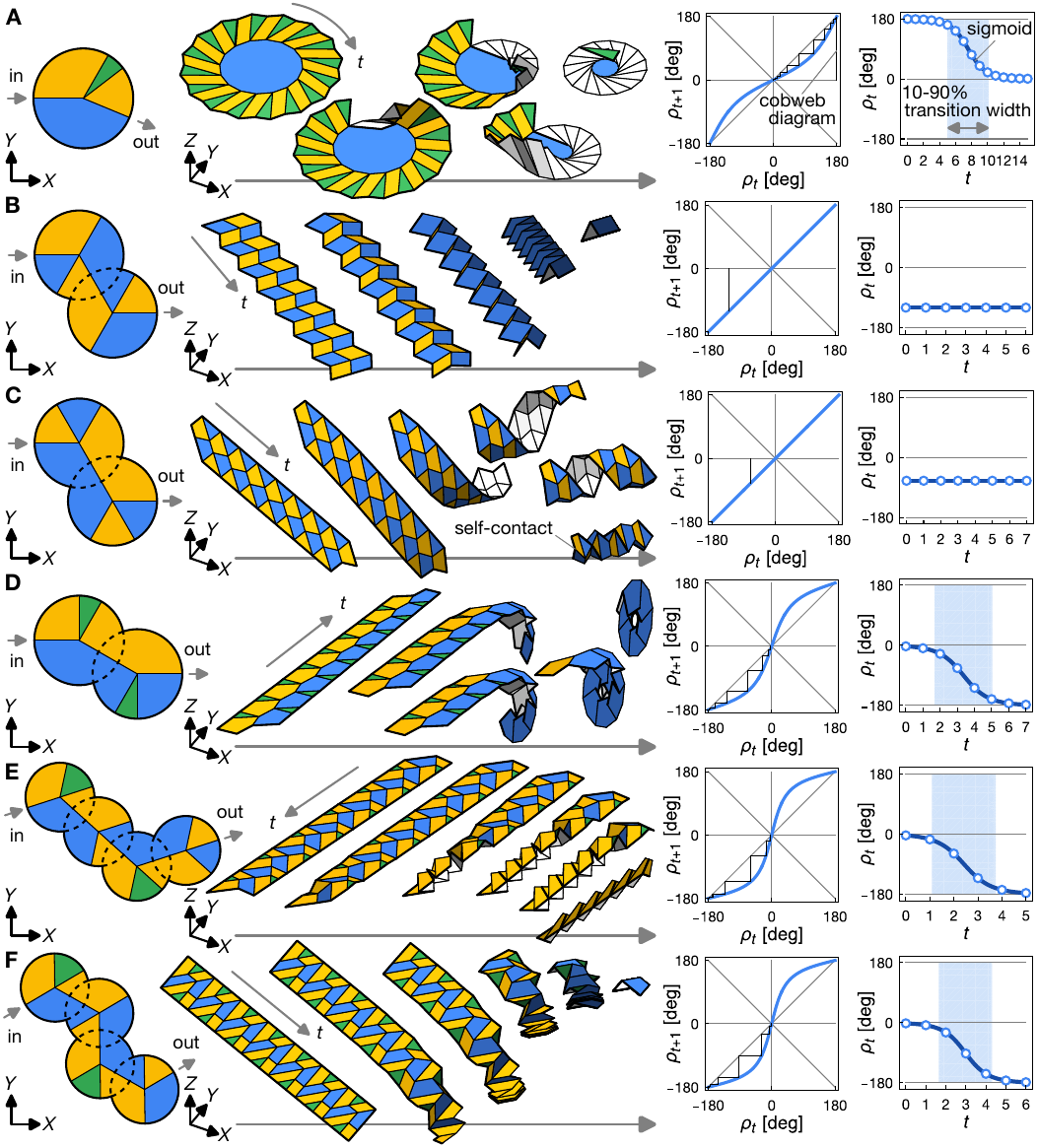}
    \caption{
        Examples of the periodic degree-$4$ strips, exhibiting the uniform or domino-like deployment.
        Each panel includes the unit-cell design, folding motion between a flat-folded state and a developed state, graph of $\rho_{t+1}=f(\rho_{t})$ with a cobweb diagram for one intermediate folded state, and the plot of the input fold angles $((t,\rho_{t}))_{t=0,1\dots}$ for that state.
        In the plots of the orbits, the underlying sigmoid function on which the trajectories lie is shown, with its 10--90\% transition width highlighted by shading.
        In visualizing folding motions, the crease lengths are specified, which do not affect the kinematic coupling behavior.
        (A) $N=1$ and $((\theta_{0}^{0},\theta_{0}^{1}),i^{\mathrm{out}}_{0}, \sigma_{0})=(
        (148.75^\circ \mathbin{,} 60^\circ)\mathbin{,} 1\mathbin{,} 1)$.
        (B) $N=2$ and $((\theta^{0}_{n},\theta^{1}_{n})\mathbin{,}i^{\mathrm{out}}_{n}\mathbin{,}\sigma_{n})_{n=0,1}=(((60^\circ,60^\circ)\mathbin{,}2\mathbin{,}1)\mathbin{,}
        ((120^\circ,120^\circ)\mathbin{,}2\mathbin{,}1))$.
        (C) $N=2$ and $((\theta^{0}_{n},\theta^{1}_{n})\mathbin{,}i^{\mathrm{out}}_{n}\mathbin{,}\sigma_{n})_{n=0,1}=(((120^\circ,120^\circ)\mathbin{,}1\mathbin{,}1)\mathbin{,}
        ((120^\circ,60^\circ)\mathbin{,}3\mathbin{,}-1))$.
        The folding motion stops before reaching the flat-folded state, because of the self-contact.
        (D) $N=2$ and $((\theta^{0}_{n},\theta^{1}_{n})\mathbin{,}i^{\mathrm{out}}_{n}\mathbin{,}\sigma_{n})_{n=0,1}=(((150^\circ,90^\circ)\mathbin{,}1\mathbin{,}-1)\mathbin{,}
        ((90^\circ,30^\circ)\mathbin{,}3\mathbin{,}-1))$.
        (E) $N=4$ and $((\theta^{0}_{n},\theta^{1}_{n})\mathbin{,}i^{\mathrm{out}}_{n}\mathbin{,}\sigma_{n})_{n=0,1,2,3}=(((120^\circ,60^\circ)\mathbin{,}1\mathbin{,}-1)\mathbin{,}
        ((120^\circ,60^\circ)\mathbin{,}2\mathbin{,}-1)\mathbin{,}
        ((120^\circ,60^\circ)\mathbin{,}3\mathbin{,}-1)\mathbin{,}
        ((120^\circ,60^\circ)\mathbin{,}2\mathbin{,}-1))$.
        (F) $N=4$ and $((\theta^{0}_{n},\theta^{1}_{n})\mathbin{,}i^{\mathrm{out}}_{n}\mathbin{,}\sigma_{n})_{n=0,1,2,3}=(((120^\circ,60^\circ)\mathbin{,}1\mathbin{,}-1)\mathbin{,}
        ((60^\circ,60^\circ)\mathbin{,}2\mathbin{,}1)\mathbin{,}
        ((120^\circ,60^\circ)\mathbin{,}3\mathbin{,}-1)\mathbin{,}
        ((120^\circ,120^\circ)\mathbin{,}2\mathbin{,}1))$.
    }
    \label{fig:examples}
\end{figure*}
For a $N$-periodic strip with $(((\theta^{0}_{n},\theta^{1}_{n}),i^{\mathrm{out}}_{n},\sigma_{n}))_{n=0,\dots,N-1}$, we define the composite map $f\coloneq f_{N-1}\circ\cdots\circ f_{0}$ which maps the input fold angle of one unit cell, $\rho_t\coloneq\rho_{Nt}^{0}$, to that of the next cell, $\rho_{t+1}$ (Fig.~\ref{fig:parameterization}(B)).
As the simplest example, when $N=1$ and $i_{0}^\mathrm{out}=1$ or $3$, the map $f$ reduces to the relation between adjacent fold angles in Eq.~\eqref{eq:fold-angle}.
The folding motion of the strip in this class is visualized in Fig.~\ref{fig:examples}(A), where the asymmetric response accumulates along the strip, resulting in domino-like deployment (see SI Appendix, Supporting Text for computational details and Movie S1 for the folding motion).
By interpreting the cell index $t$ as discrete time, the map $f$ defines a discrete dynamical system in which the developed state ($\rho=0$) and the flat-folded state ($\rho=\pi$ or $-\pi$) are fixed points.
The observed domino-like deployment corresponds to heteroclinic orbits connecting these fixed points, i.e., trajectories that approach one fixed point as $t \to -\infty$ and the other as $t \to +\infty$.
Notably, Eq.~\eqref{eq:fold-angle} implies that the cosines of adjacent fold angles are related by a fractional linear transformation.
Such transformations are characterized, up to an overall scaling, by the ratio $A/B$, which therefore serves as the effective parameter of the map.
Furthermore, the orbit $(\rho_{t})_{t=0,\dots}$ lies on the function:
\begin{align}
    \begin{split}
        &2\arctan{\left(\left(\tan{\dfrac{\rho_{0}}{2}}\right)p^{t}\right)},\\
        \quad\text{where}\quad p&:=-\mathrm{sgn}(\cos{\theta^{0}+\sigma\cos{\theta^{1}}})\sqrt{\frac{A-B}{A+B}}.
    \end{split}
\end{align}
This function is a continuous sigmoid on $\mathbb{R}$ for $p>0$; for $p<0$, the same property holds after taking its absolute value.
By examining the 10--90\% transition width of this sigmoid function, $|2\log(\tan(\pi/20))/\log{p}|$, one can estimate how many unit cells are required for the domain wall connecting the developed and flat-folded states.
For example, in Fig.~\ref{fig:examples}(A), the transition width corresponds to approximately five unit cells.
Although the above concerns the case of $N=1$, the class of symmetric fractional linear transformations is closed under composition, so that the composite map $f$ retains the same functional form; thus, the same argument applies to $N>1$.
Moreover, provided that the map is not the identity (i.e., $B\neq 0$), the fixed points are universally given by $\cos\rho^0=\pm1$, independent of the parameters, while their stability, as quantified by the slope at the fixed points, depends on the parameters (see also Fig.~\ref{fig:asymmetry}(B)).
Such degenerate cases arise, for example, when $i_n^{\mathrm{out}}=2$ for all $n$ (Fig.~\ref{fig:examples}(B); Movie S2), or when the asymmetry is canceled by specific choices of sector angles (Fig.~\ref{fig:examples}(C); Movie S3).
In many developable and flat-foldable quadrilateral-mesh origami, such as the Miura-ori, the sub-strips fall into these degenerate cases; consequently, most origami-based mechanisms based on these patterns exhibit uniform deformation.
However, except for these singular cases, domino-like deployment generically occurs when adjacent creases are used.
Fig.~\ref{fig:examples}(D)--(F) present representative generic examples without the special symmetries of the degenerate cases (Movies S4--S6).
Fig.~\ref{fig:examples}(D) has the same connectivity as Fig.~\ref{fig:examples}(C), but the sector angles are no longer symmetric.
Fig.~\ref{fig:examples}(E) and (F) correspond to designs in which vertices are inserted into a sub-strip of a Miura-ori-like pattern via adjacent-vertex connections.
In particular, Fig.~\ref{fig:examples}(F) exhibits the folding sequence similar to that of Fig.~\ref{fig:examples}(C), suggesting the high packaging ratio.
In summary, periodic degree-4 origami strips generically exhibit sequential deployment under adjacent-crease connections (except for degenerate cases), with the transition width---and hence the number of involved unit cells---controlled by the effective parameter $A/B$ of the composite map $f$.

\subsection*{Programmable macroscopic shape while preserving propagation behavior}
\begin{figure*}[htbp]
    \centering
    \includegraphics[width=\linewidth]{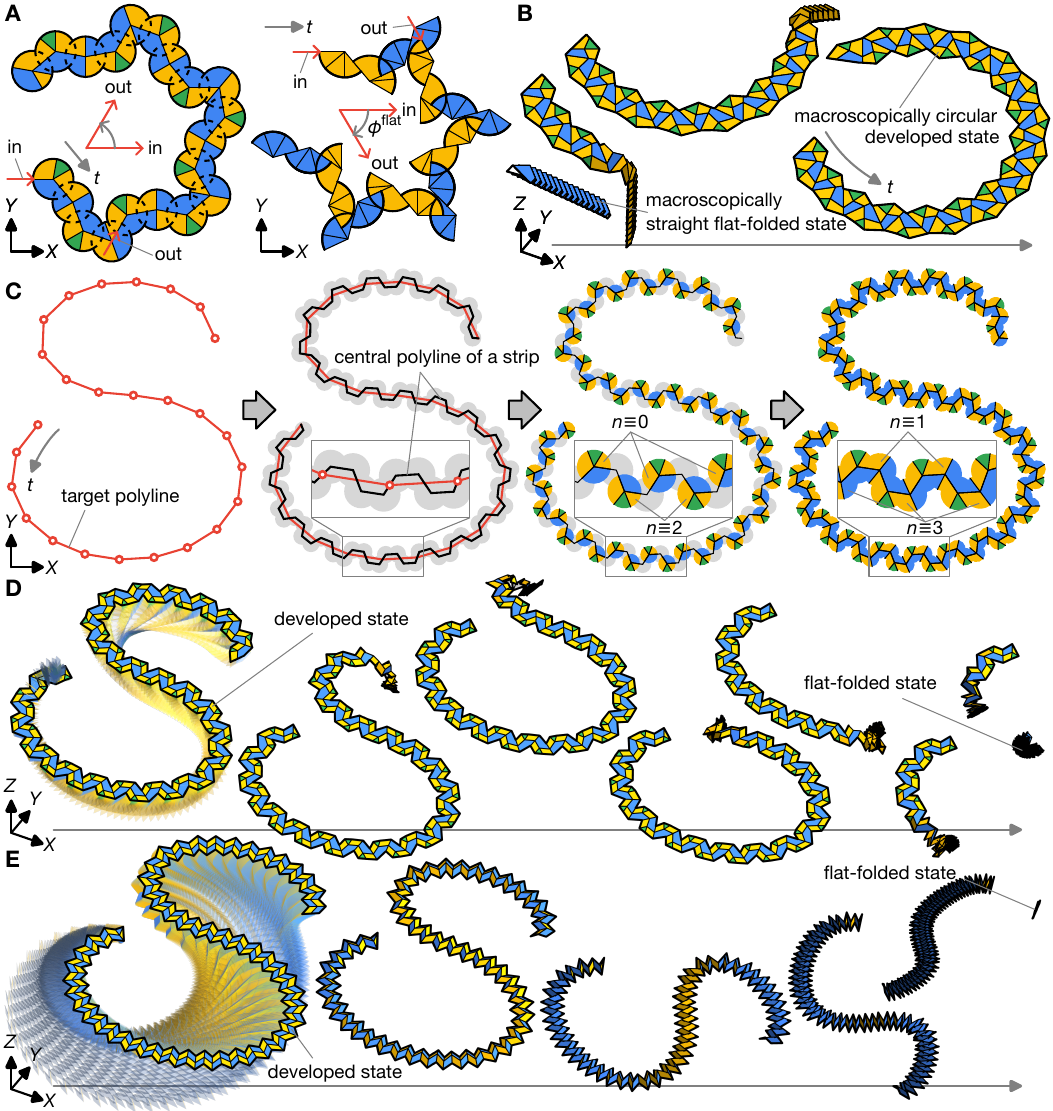}
    \caption{
    (A) Illustration of the total turning angle over one unit cell ($N=4$) in the developed state (left) and the flat-folded state (right).
    (B) Domino-like deployment of the periodic strip from the straight flat-folded state to the circular developed state ($N=4\mathbin{,}((\theta_{n}^{0},\theta_{n}^{1})\mathbin{,}i_{n}^{\mathrm{out}}\mathbin{,}\sigma_{n})_{n=0,1,2,3}=(((110^\circ,70^\circ)\mathbin{,}1\mathbin{,}-1)\mathbin{,}((60^\circ,60^\circ)\mathbin{,}2\mathbin{,}1)\mathbin{,}((110^\circ,70^\circ)\mathbin{,}3\mathbin{,}-1)\mathbin{,}((130^\circ,130^\circ)\mathbin{,}2\mathbin{,}1))$).
    (C) Design procedure for a non-periodic strip that deploys to a prescribed target shape.
    The period, connectivity, mode assignment, and discrete dynamical system are identical to those in Fig.~\ref{fig:examples}(F).
    Given a target deployed shape represented as a polyline, each segment is first mapped to the central polyline of a unit cell.
    The sector angles of the vertices with $n \equiv 0,2$ are then determined from $\theta_0^0$ and $\theta_2^1$, which are fixed by the central polyline, together with the prescribed effective parameter $A/B$.
    The sector angles of the vertices with $n \equiv 1,3$ are subsequently determined.
    Although each vertex admits one degree of freedom in its sector angles, mirror-symmetric ones are selected.
    (D), (E) Folding motions of strips that deploy into an S-shape via domino-like (top) and uniform (bottom) deployment. 
    To compare swept volumes, snapshots of the folding motion are overlaid on the developed configurations.
    }
    \label{fig:curvature}
\end{figure*}

The developed and flat-folded states of the examples in Fig.~\ref{fig:examples}(E) and (F) are both macroscopically straight when viewed along the central polyline connecting the vertices.
More generally, for periodic strips, the macroscopic shapes of the developed and flat-folded states are characterized by the total turning angle over one unit cell (Fig.~\ref{fig:curvature}(A)),
\begin{align}\label{eq:curvature}
    \begin{split}
        \phi^{\mathrm{dev}}&\coloneq\sum_{n=0}^{N-1}\left(\pi+\sum_{i=1}^{i_n^{\mathrm{out}}}\theta_{n}^{i-1}\right)\mod{2\pi},\\
    \quad \phi^{\mathrm{flat}}&\coloneq\sum_{n=0}^{N-1}\left(\pi+\sum_{i=1}^{i_n^{\mathrm{out}}}(-1)^{s(n,i)}\theta_{n}^{i-1}\right)\mod{2\pi},
    \end{split}\notag
\end{align}
\begin{align}
    \begin{split}
        \text{where}\quad &s(n,i)\coloneq i+\sum_{m=0}^{n-1}\left(i_{m}^{\mathrm{out}}-1\right)\mod{2},\\
    &\phi^{\mathrm{dev}},\phi^{\mathrm{flat}}\in(-\pi,\pi].
    \end{split}
\end{align}
The exponent $s(n, i)\equiv 0$ ($s(n, i)\equiv 1$) indicates that the $i$-th face of the $n$-th vertex is face-up (face-down) at the flat-folded state.
The quantities $\phi^{\mathrm{dev}}$ and $\phi^{\mathrm{flat}}$ depend only on the connectivities and sector angles and are independent of the crease lengths between the vertices, i.e., they represent a discrete curvature accumulated over one unit cell.
The examples in Fig.~\ref{fig:examples}(E) and (F) correspond to the case $\phi^{\mathrm{dev}}=\phi^{\mathrm{flat}}=0$.
By tuning sector angles to control $\phi^{\mathrm{dev}}$ and $\phi^{\mathrm{flat}}$, one can realize mechanisms that deploy in a domino-like manner, for example, from a straight flat-folded state into a circular developed state (Fig.~\ref{fig:curvature}(B)).
In periodic strips, this total turning angle is constant across cells, yielding straight or circular configurations.
To realize more general shapes, it is necessary to introduce spatial variations in the sector angles.
However, such variations typically render the recurrence relation non-autonomous, i.e., dependent on the cell index $t$, leading to spatially varying propagation behavior.
Importantly, as shown in Eq.~\eqref{eq:fold-angle}, the relation between opposite fold angles reduces to the identity; therefore, variations in the sector angles at vertices with $i_n^{\mathrm{out}}=2$---for example, the $n\equiv1,3$ vertices in Fig.~\ref{fig:examples}(E) and (F)---do not affect the discrete dynamical system $f$. 
Consequently, by using these sector angles as design parameters, one can tune the macroscopic curvature independently while preserving the propagation characteristics.
Furthermore, even for vertices with $i_n^{\mathrm{out}}=1$ or $3$, the recurrence relation can be preserved under variations of the sector angles that keep the effective parameter $A/B$ invariant.
This principle extends to non-periodic designs in which only the sector angles are varied while the connectivity and mode-assignment are kept fixed. 
Specifically, we can make the sector angles through the strip non-periodic, while the recurrence relation $f$ remains autonomous, i.e., independent of the cell index $t$.
As a result, one can design strips that follow arbitrary planar curves in their developed state while preserving the propagation behavior (Fig.~\ref{fig:curvature}(C) and (D); Movie S7; see SI Appendix, Supporting Text for details of mapping the target polyline to the central crease).
Compared with a strip based on uniform deployment (Fig.~\ref{fig:curvature}(E); Movie S8), the present design exhibits a more localized deformation during deployment, suggesting a reduced swept volume.
These results demonstrate that the macroscopic shape of domino-deploying origami strips can be programmed independently of their propagation behavior.

\subsection*{Demonstration of the domino-like deployment}
\begin{figure*}[htbp]
    \centering
    \includegraphics[width=\linewidth]{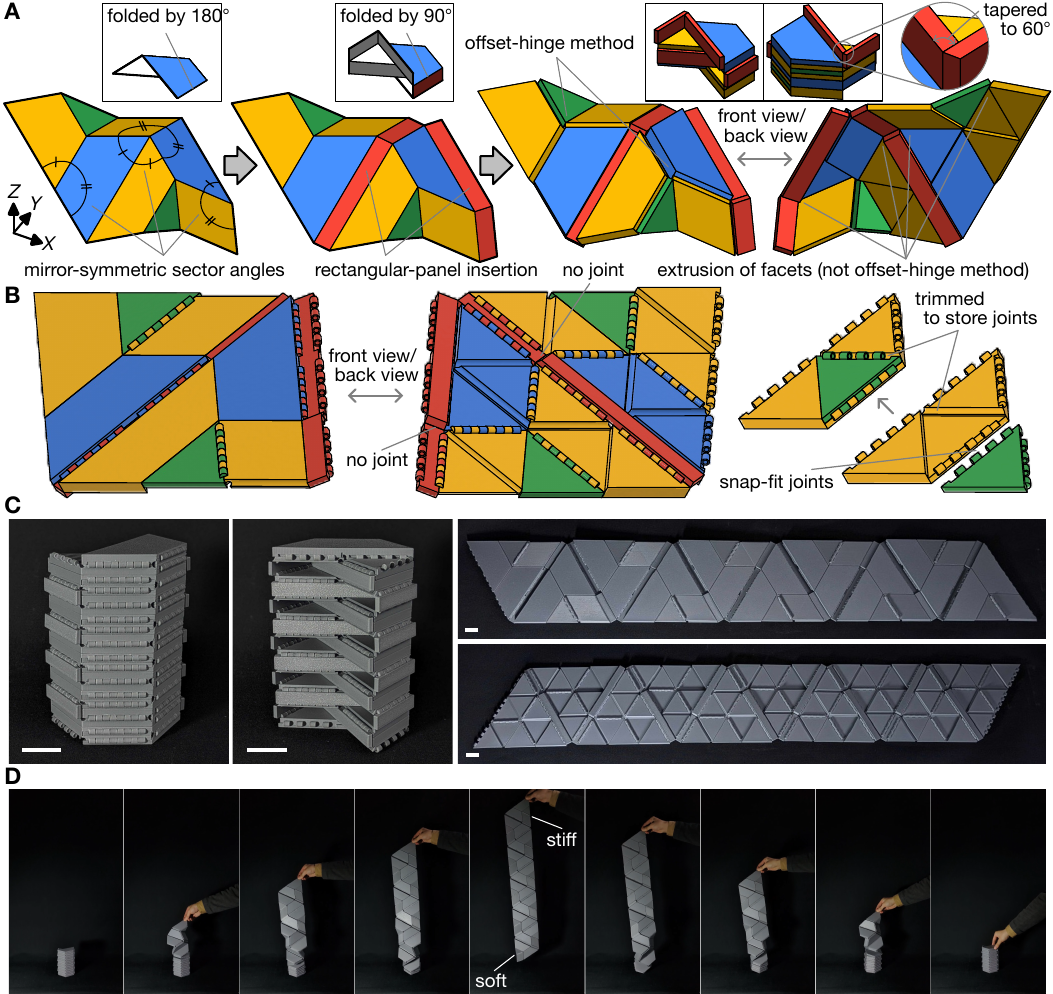}
    \caption{
        (A) Transformation process of a partially folded zero-thickness unit cell via rectangular-panel insertion and panel thickening using the face-extrusion and offset-hinge methods.
        Insets show the flat-folded states at each step.
        (B) Unit cell assembled from panels with snap-fit hinges placed along the creases.
        For simplicity, joints at the short creases between rectangular panels are omitted, as they do not affect the kinematics.
        Some panel corners are trimmed to accommodate the joints in the flat-folded state.
        (C) Prototype of the strip consisting of five unit cells (left: front and bottom view of the flat-folded state; right: front and bottom view of the developed state).
        Scale bar indicates 2cm.
        Each panel was 3D printed using a Bambu Lab P1S with PLA filament.
        (D) Snapshots of the deployment and stowage of the prototype.
    }
    \label{fig:demonstration}
\end{figure*}
To demonstrate the feasibility of domino-like deployment, we constructed a physical prototype based on the design shown in Fig.~\ref{fig:examples}(F).
To this end, we transformed the zero-thickness origami mechanism into a structure composed of panels with finite thickness while preserving its kinematics.
Accounting for panel thickness prevents unintended elastic deformation of facets that may interfere with the kinematic coupling between vertices---an effect beyond the scope of this study---and is also essential for practical engineering applications such as deployable structures and robotics.
Here, we employ the offset-hinge method, which is applicable to developable and flat-foldable degree-4 vertices~\cite{hoberman1988reversibly,chen2015origami,lang2018review,zhang2018mobile}.
In this approach, each panel is assigned a thickness determined by the sector angles $(\theta^{0},\theta^{1})$, and hinges are placed on the bottom (top) surfaces for mountain (valley) folds.
A straightforward application of this method to each vertex in the strip shows that avoiding self-intersection requires the panel thickness to vary exponentially along the strip, rendering the design impractical (see SI Appendix, Supporting Text for details).
To resolve this issue, we introduce rectangular panels at two of the four vertices in each unit cell (Fig.~\ref{fig:demonstration}(A)), which do not change the kinematics because of the mirror-symmetry of sector angles.
At the vertices where rectangular panels are introduced, the corresponding creases are effectively split into two parallel creases, and they are folded at $90^\circ$ in the flat-folded state.
This feature allows these vertices to be thickened without relying on the offset-hinge method; instead, the facets can be simply extruded along their normals.
Consequently, the offset-hinge method needs to be applied only to the remaining two vertices.
This strategy yields a design with uniform panel thickness along the strip, although the deployed configuration becomes staircase-like (Movie S9).
Finally, snap-fit hinges are incorporated into each panel, and local trimming is applied where necessary to prevent interference during folding (Fig.~\ref{fig:demonstration}(B)).
The panels are fabricated by 3D printing and assembled into a prototype (Fig.~\ref{fig:demonstration}(C)).
The structure exhibits the expected behavior: when deployed from one end, the transition propagates sequentially along the strip, reaching the fully deployed state (Fig.~\ref{fig:demonstration}(D); Movie S10).
Conversely, during folding, the structure returns to the flat-folded state in the reverse order of deployment.
Moreover, the opposite end remains mechanically stiff at the developed state, making it difficult to initiate folding from that side, which is analogous to topological polarization in mechanical systems~\cite{chen2014nonlinear,chen2016topological,rocklin2017transformable,mao2018maxwell}.
In zero-thickness origami, multiple folding branches can coexist at the developed state, requiring an explicit assignment of the folding mode, which is specified by $\sigma$ in our parameterization.
In contrast, introducing panel thickness imposes geometric constraints that naturally suppress some of these modes: owing to the offset-hinge design, each crease is effectively biased to fold either as a mountain or a valley, thereby reducing the admissible configurations.
In the present design, a residual mode remains at the developed configuration in which folding occurs along collinear creases of the rectangular panels.
Nevertheless, we confirm that, once folding is initiated from one end in the intended manner, the structure sequentially follows the desired mode throughout the strip.
Although we fabricated a prototype for a specific design, the presented thickness-accommodation method extends to other patterns. 
In particular, mirror symmetry of the sector angles at vertices with $n\equiv 1,3$ is a necessary and sufficient condition for inserting rectangular panels without altering the kinematics. 
Under this condition, the strips in Fig.~\ref{fig:curvature} (B) and (D) can be thickened in the same manner.
In practice, however, the required panel trimming and hinge placement may become more intricate depending on the geometry.
Despite these practical complexities, the qualitative agreement between the designed kinematics and the observed behavior demonstrates that domino-like deployment can be realized in physical origami structures.

\section*{Discussion}
In this work, we established a design framework for kinematic propagation of transition fronts in geometrically constrained systems, based on their correspondence with heteroclinic orbits in discrete dynamical systems.
Using strips of developable, flat-foldable degree-4 origami vertices, we showed that asymmetric kinematic coupling between adjacent creases gives rise to domino-like sequential deployment without energy barriers, and that macroscopic curvature can be programmed independently of the propagation behavior.
In addition, the localized nature of the propagation suggests that such mechanisms may be advantageous for deployment under spatial constraints, where minimizing swept volume is critical.
We further validated this concept through a thick-panel origami prototype exhibiting the designed deployment.

Relaxing the developability condition may enable deployment from compact, flat-folded configurations to spatial shapes with nonzero torsion, and extending the architecture to higher-dimensional structures, such as sheets and cellular assemblies, remains an important direction for future work.
Such extensions may, however, require additional design considerations, as many origami-based mesh or cellular structures (e.g., Miura-ori) are overconstrained and tend to exhibit uniform rather than domino-like deformation.

Importantly, because the propagation behavior is governed by an underlying discrete dynamical system, geometric parameters can be modified while preserving this structure.
This enables design flexibility beyond macroscopic curvature, such as tailoring point trajectories during deployment, minimizing swept volume, or tuning kinematic responses to specific functional requirements, without altering the propagation characteristics.

More broadly, linking local kinematic coupling to global propagation through discrete dynamical systems provides a general framework for understanding and designing transition fronts in geometrically constrained systems.
By advancing from isolated examples toward a systematic design methodology, this work opens new possibilities for deployable structures, robotics, and mechanical metamaterials.

\begin{acknowledgments}
    R.I. acknowledges funding from JSPS KAKENHI Grant No. JP26KJ0431.
    R.I. and T.T. acknowledge funding from JSPS KAKENHI Grant No. JP24H00822.
\end{acknowledgments}

\bibliography{myref}

\clearpage
\onecolumngrid
\appendix
\input{appendix.tex}

\end{document}

%% file: appendix.tex
\section{Computation and Visualization of Folded Configurations}
We describe how the folded configurations of the strips shown in Figs.~3 and 4 in the main text are visualized in three-dimensional space.
The folded state of a strip is computed recursively, in a manner analogous to the fold angles, by sequentially determining the configuration of each degree-4 vertex.
Although the strips considered in the main text are developable and flat-foldable, the representation and computational procedure described below are applicable more generally to vertices and strips that do not necessarily satisfy these conditions.
Accordingly, we first describe the representation and computation of the folded state of a single degree-4 vertex, followed by those for the entire strip.

\subsection{Vertex configuration}
The folded state of a degree-4 vertex with sector angles $(\theta^{i})_{i=0,1,2,3}$ can be determined by specifying the fold angle $\rho^{0}$ and computing the remaining fold angles $(\rho^{i})_{i=1,2,3}$ according to Eq.~(1) in the main text.
However, to specify its spatial configuration in three-dimensional space, both position and orientation must be defined.
Here, we represent the position and orientation using the following quantities: the position of the central vertex $\mathbf{o} \in \mathbb{R}^3$, unit vectors $\mathbf{c}^{i} \in \mathbb{R}^3$ indicating the directions of the creases, and unit normal vectors $\mathbf{n}^{i} \in \mathbb{R}^3$ of the faces.
Each vector $\mathbf{c}^{i}$ is based at $\mathbf{p}$.
The vector $\mathbf{n}^{i}$ denotes the unit normal of the face bounded by the $i$-th and $(i+1)$-th creases, and the configuration of this face is represented by the triplet $(\mathbf{o}, \mathbf{c}^{i}, \mathbf{n}^{i})$.
Once the folded state is specified, i.e., when $(\rho^{i})_{i=0,1,2,3}$ are given, the configurations of the four faces can be computed recursively from the configuration of the $0$-th face, $(\mathbf{o}, \mathbf{c}^{0}, \mathbf{n}^{0})$, using the relations
\begin{equation}\label{eq:configuration}
	\mathbf{c}^{i+1}=\mathbf{R}(\theta^{i},\mathbf{n}^{i})\mathbf{c}^{i}\quad \text{and}\quad\mathbf{n}^{i+1}=\mathbf{R}(\rho^{i+1},\mathbf{c}^{i+1})\mathbf{n}^{i},
\end{equation}
where $\mathbf{R}(\theta, \mathbf{v}) \in \mathbb{R}^{3 \times 3}$ denotes the rotation matrix representing a rotation by angle $\theta$ about the axis $\mathbf{v}$.
From the above, the folded state and spatial configuration of a degree-4 vertex can be fully specified by $\rho^{0}$ together with $(\mathbf{o}, \mathbf{c}^{0}, \mathbf{n}^{0})$.

\subsection{Strip configuration}
Next, we consider the configuration of a strip specified by the sequence of sector angles and output-crease indices, $((\theta_{n}^{i})_{i=0,1,2,3}\mathbin{,}i_{n}^{\mathrm{out}})_{n=0,1,\dots}$.
We also assign the crease lengths at the $n$-th vertex as $(l_{n}^{i})_{i=0,1,2,3}$ where $l_{n}^{i}\in\mathbb{R}_{>0}$; from the connectivity between adjacent vertices, these satisfy $l_{n+1}^{0} = l_{n}^{i_{0}^\mathrm{out}}$.
The folded state of the strip is determined by specifying the input fold angle $\rho_{0}^{0}$ at the $0$-th vertex and computing all other fold angles using Eq.~(1) in the main text.
Following the single-vertex formulation, we represent the configuration of the $n$-th vertex by $\mathbf{o}_{n}$, $(\mathbf{c}_{n}^{i})_{i=0,1,2,3}$, and $(\mathbf{n}_{n}^{i})_{i=0,1,2,3}$, and use these sequences to describe the configuration of the entire strip.
The configuration of the $0$-th vertex can be computed from Eq.~\eqref{eq:configuration} once $\mathbf{c}_{0}^{0}$ and $\mathbf{n}_{0}^{0}$ are specified.
From the connectivity between the $0$-th and $1$-st vertices, we have
\begin{equation}
    \mathbf{o}_{1} = \mathbf{o}_{0} + l_{0}^{i_{0}^{\mathrm{out}}}\,\mathbf{c}_{0}^{i_{0}^{\mathrm{out}}},\quad\text{and}\quad
    \mathbf{c}_{1}^{0} = -\mathbf{c}_{0}^{i_{0}^{\mathrm{out}}}.
\end{equation}
Furthermore, since the $0$-th face of the $1$-st vertex lies in the same plane as the $(i_{0}^{\mathrm{out}}-1)$-th face of the $0$-th vertex, we have
\begin{equation}
    \mathbf{n}_{1}^{0} = \mathbf{n}_{0}^{i_{0}^{\mathrm{out}}-1}.
\end{equation}
Thus, the configuration of the $1$-st vertex can be computed from Eq.~\eqref{eq:configuration}.
Generalizing this procedure, the parameters representing the configuration of the $(n+1)$-th vertex are determined from those of the $n$-th vertex as
\begin{equation}\label{eq:configuration_recurrence}
    \mathbf{o}_{n+1}=\mathbf{o}_{n}+l_{n}^{i_{n}^\mathrm{out}}\mathbf{c}_{n}^{i_{n}^{\mathrm{out}}},\quad
    \mathbf{c}_{n+1}^{0}=-\mathbf{c}_{n}^{i_{n}^{\mathrm{out}}},\quad\text{and}\quad
    \mathbf{n}_{n+1}^{0}=\mathbf{n}_{n}^{i_{n}^{\mathrm{out}}-1}.
\end{equation}
By this recursive procedure, the folded state and spatial configuration of the strip can be fully determined from $\rho_{0}^{0}$ and $(\mathbf{o}_{0}, \mathbf{c}_{0}^{0}, \mathbf{n}_{0}^{0})$, using Eq.~\eqref{eq:configuration} and Eq.~\eqref{eq:configuration_recurrence}.
In Figs.~3 and 4 of the main text, each strip is visualized by placing points along the lateral creases according to the crease lengths $(l_{n}^{i})_{i=0,1,2,3}$, and constructing faces by connecting sets of coplanar points.

\section{Thickness Accommodation of the Strip by Offset-hinge Method}
Here, we first briefly introduce the offset-hinge method~\cite{hoberman1988reversibly,chen2015origami,zhang2018mobile,lang2018review}, a thickness-accommodation technique for developable and flat-foldable degree-4 vertices, focusing on its underlying principle and design conditions.
We then show that, when the offset-hinge method is applied to each vertex of the strip shown in Fig.~3(F) in the main text, the panel thickness varies exponentially along the strip.

\subsection{Introduction to offset-hinge method}
Following previous work~\cite{chen2015origami}, we introduce the offset-hinge method.
Fig.~\ref{fig:offset-hinge}(A) and (B) show an example of applying the offset-hinge method to a degree-4 vertex.
Focusing on the relative heights of the hinges in the developed state, the valley crease lies at the highest position, the mountain crease opposite to it lies at the lowest position, and the remaining mountain creases are located at intermediate heights.
In the figure, the distances between adjacent creases, $(d^{i})_{i=0,1,2,3}$, determine the thickness of the panels.
Importantly, the values of $(d^{i})_{i=0,1,2,3}$ are not arbitrary; to preserve kinematic equivalence, they must satisfy
\begin{equation}\label{eq:linklength}
d^{1} = d^{0}\frac{\sin{\theta^{1}}}{\sin{\theta^{0}}}, \quad
d^{2} = d^{0}, \quad\text{and}\quad
d^{3} = d^{1}.
\end{equation}
Eq.~\eqref{eq:linklength} arises from the correspondence between origami kinematics and linkage mechanisms: zero-thickness origami, in which creases intersect at a point, can be regarded as a 4R spherical linkage, whereas thickened origami, in which creases no longer intersect at a point, corresponds to a 4R spatial linkage.
In this correspondence, the design parameters $(\theta^{i})_{i=0,1,2,3}$ and $(d^{i})_{i=0,1,2,3}$ correspond to the twist angles and link lengths of the 4R spatial linkage, respectively.
Unlike a 4R spherical linkage, which is mobile for arbitrary sector angles, a 4R spatial linkage is generally overconstrained and exhibits mobility only when its twist angles and link lengths satisfy specific conditions.
The developability and flat-foldability conditions of a single-vertex origami, together with Eq.~\eqref{eq:linklength}, correspond to one such condition under which the 4R spatial linkage becomes an overconstrained but movable mechanism, known as the Bennett linkage.
Therefore, when thickening a developable and flat-foldable degree-4 vertex using the offset-hinge method, the design has only one independent parameter, namely a single link length.

\subsection{Exponential change of panel-thickness along the strip}
\begin{figure}[t]
    \centering
    \includegraphics[width=\textwidth]{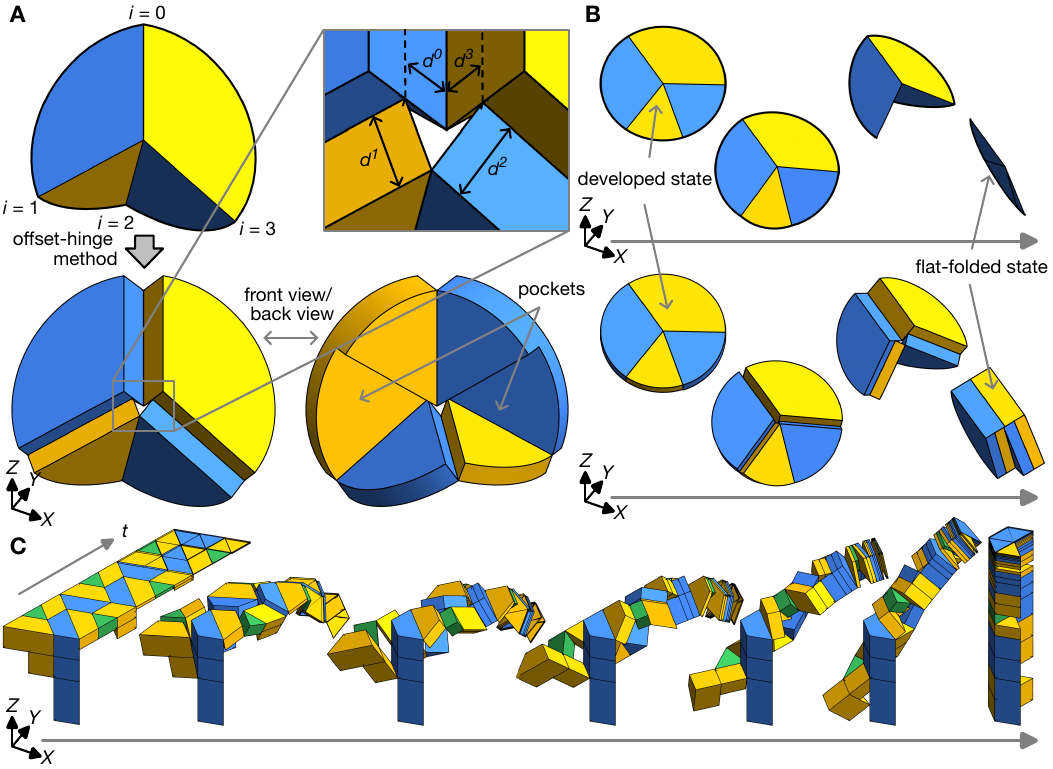}
    \caption{
        (A) Thickening of a developable and flat-foldable degree-4 vertex using the offset-hinge method.
        The panel thickness is determined by the offset distances between adjacent creases, parameterized by $(d^{i})_{i=0,1,2,3}$.
        (B) Equivalence of the folding motions before and after applying the offset-hinge method.
        (C) Snapshots of the folding motion of the strip exhibiting an exponential variation in panel thickness.
        Despite consisting of only three unit cells, the panel thickness varies significantly along the strip.
    }
    \label{fig:offset-hinge}
\end{figure}

When the offset-hinge method is applied to each vertex of the strip shown in Fig.~3(F) in the main text, the representative link lengths at each vertex, $(l_{n}^{0})_{n=0,1,\dots}$, serve as the design parameters.
In structures such as quadrilateral meshes, where internal vertices form loops, the link lengths must satisfy certain compatibility conditions in order to preserve mobility~\cite{zhang2018mobile}.
In contrast, a strip contains no such loops; therefore, mobility is preserved for arbitrary choices of $(l_{n}^{0})_{n=0,1,\dots}$.
However, the situation changes when self-intersection between panels is taken into account.
To avoid self-intersection in the flat-folded state, the developed configuration must be geometrically flat.
As a result, once $l_{0}^{0}$ is specified, all other link lengths are uniquely determined.
Consequently, the panel thickness is governed by the relative heights of the hinges in the deployed state.
For the connectivity and mountain–valley assignment considered here, the relative height between the input and output creases remains unchanged at vertices with $n \equiv 1,3$, where opposite creases connect adjacent vertices.
In contrast, at vertices with $n \equiv 0,2$, where adjacent creases are used for connectivity, the output crease is always located at a higher position than the input crease.
As a result, the hinge positions increase monotonically with the cell index.
The combination of these two effects leads to an exponential decrease in panel thickness along the strip as the cell index increases (Fig.~\ref{fig:offset-hinge}(C)).
The thickness-accommodation method presented in the main text avoids this exponential variation by employing extrusion, which eliminates the requirement that the deployed configuration be geometrically flat.

\section{Mapping a Target Polyline to the Central Crease of the Strip}
\begin{figure}[t]
    \centering
    \includegraphics[width=\textwidth]{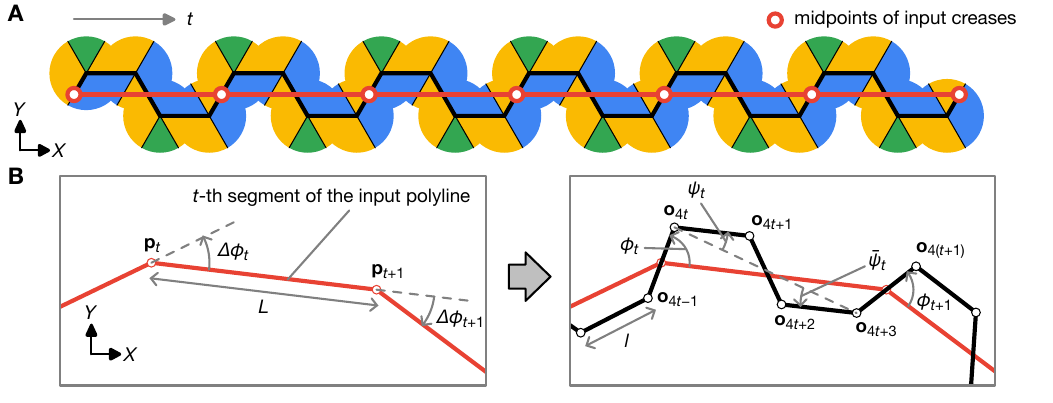}
    \caption{
        (A) Straight polyline drawn on the crease pattern of the strip in Fig.~3(F) of the main text, connecting the midpoints of the input creases.
        (B) Mapping of the $t$-th segment of the input polyline (left) into the central creases of the $t$-th cell (right).
    }
    \label{fig:mapping-polyline}
\end{figure}

We describe how the input polyline $(\mathbf{p}_{t})_{t=0,\dots,T}$ with $\mathbf{p}_{t} \in \mathbb{R}^{2}$ is transformed into the central crease $(\mathbf{o}_{n})_{n=0,\dots,4T-1}$ with $\mathbf{o}_{n} \in \mathbb{R}^{2}$ in Fig.~4(C) of the main text.
For simplicity, we assume that both the input polyline and the central crease have uniform segment lengths, denoted by $L\in\mathbb{R}_{>0}$ and $l\in\mathbb{R}_{>0}$, respectively.
As a natural construction, the $t$-th segment of the polyline is associated with the central crease of one unit cell, represented by $(\mathbf{o}_{4t}, \mathbf{o}_{4t+1}, \mathbf{o}_{4t+2}, \mathbf{o}_{4t+3})$.
We further impose, as a design requirement, that when a straight polyline is given as input, the developed configuration shown in Fig.~3(F) of the main text is reproduced (up to a suitable adjustment).
As illustrated in Fig.~\ref{fig:mapping-polyline}(A), connecting the midpoints of the input creases of each unit cell in the developed configuration of Fig.~3(F) yields a straight polyline.
Based on this observation, we assign the endpoints of each segment to the midpoints of the input and output creases of each unit cell.
Under this construction, the mapping of the $t$-th segment to the central crease of the $t$-th cell can be parameterized by two quantities:
(1) the turning angle $\phi_{t}$ from the unit direction vector of the $t$-th segment, $\mathbf{v}_{t}:=(\mathbf{p}_{t+1}-\mathbf{p}_{t})/\|\mathbf{p}_{t+1}-\mathbf{p}_{t}\|$, to the vector $(\mathbf{o}_{4t} - \mathbf{p}_{t}) / \|\mathbf{o}_{4t} - \mathbf{p}_{t}\|$, and
(2) the turning angle $\psi_{t}$ from $(\mathbf{o}_{4t+3} - \mathbf{o}_{4t})/\|\mathbf{o}_{4t+3} - \mathbf{o}_{4t}\|$ to $(\mathbf{o}_{4t+1} - \mathbf{o}_{4t}) / \|\mathbf{o}_{4t+1} - \mathbf{o}_{4t}\|$ (Fig.~\ref{fig:mapping-polyline}(B)).
Using these parameters, the positions of the four vertex centers within the $t$-th unit cell in the developed state are given by
\begin{align}\label{eq:four-points}
    \begin{split}
        \mathbf{o}_{4t} &= \mathbf{p}_{t} + \frac{l}{2}\,\mathbf{R}(\phi_{t})\,\mathbf{v}_{t},\\
        \mathbf{o}_{4t+3} &= \mathbf{p}_{t+1} - \frac{l}{2}\,\mathbf{R}(\phi_{t+1})\,\mathbf{v}_{t+1},\\
        \mathbf{o}_{4t+1} &= \mathbf{o}_{4t} + l\,\mathbf{R}(\psi_{t})\left(\frac{\mathbf{o}_{4t+3} - \mathbf{o}_{4t}}{\|\mathbf{o}_{4t+3} - \mathbf{o}_{4t}\|}\right),\\
        \mathbf{o}_{4t+2} &= \mathbf{o}_{4t+3} + l\,\mathbf{R}(\bar{\psi}_{t})\left(\frac{\mathbf{o}_{4t} - \mathbf{o}_{4t+3}}{\|\mathbf{o}_{4t} - \mathbf{o}_{4t+3}\|}\right).
    \end{split}
\end{align}
Here, $\mathbf{R}(\phi)\in\mathbb{R}^{2\times 2}$ denotes the two-dimensional rotation matrix with rotation angle $\phi$, and $\bar{\psi}_{t}$ is defined as
\begin{align}
    \begin{split}
	\bar{\psi}_{t}&\coloneqq\arccos{\left(\left(\dfrac{a^2}{2la}\right)
    \left(\dfrac{a^2+\|\mathbf{o}_{4t+3}-\mathbf{o}_{4t}\|^2-l^2}{2a\|\mathbf{o}_{4t+3}-\mathbf{o}_{4t}\|}\right)
    +\sqrt{1-\left(\dfrac{a^2}{2la}\right)^2}\sqrt{1-\left(\dfrac{a^2+\|\mathbf{o}_{4t+3}-\mathbf{o}_{4t}\|^2-l^2}{2a\|\mathbf{o}_{4t+3}-\mathbf{o}_{4t}\|}\right)^2}\right)},\\
    \text{where}&\quad a:=\sqrt{l^2+\|\mathbf{o}_{4t+3}-\mathbf{o}_{4t}\|^2-2l\|\mathbf{o}_{4t+3}-\mathbf{o}_{4t}\|\cos{\psi_{t}}}.
    \end{split}
\end{align}
Instead of introducing and optimizing an explicit objective function, we determine $\phi_{t}$ and $\psi_{t}$ heuristically.
Specifically, $\phi_{t}$ is prescribed using the turning angle of the polyline segments as
\begin{align}\label{eq:turning-angle-1}
    \begin{split}
        \phi_{t}&=\phi^{*}+\dfrac{\Delta\phi(\mathbf{v}_{t-1},\mathbf{v}_{t})}{2},\\
        \text{where}\quad\Delta\phi([x_{1},y_{1}]^T,&[x_{2},y_{2}]^T) \coloneqq \operatorname{atan2}\!\big(x_{1} y_2 - y_1 x_2,\; x_1 x_2 + y_1 y_2 \big).
    \end{split}
\end{align}
Eq.~\eqref{eq:turning-angle-1} sets the target angle $\phi^{*}$ when the polyline is locally straight, and adjusts $\phi_{t}$ according to the deviation from straightness, measured by $\Delta\phi$.
The initial rotation angle $\phi_{0}$ at the starting point can be chosen arbitrarily.
Next, $\psi_{t}$ is prescribed by enforcing symmetry, i.e., $\psi_{t} = \bar{\psi}_{t}$, as
\begin{align}\label{eq:turning-angle-2}
    \begin{split}
    	\psi_{t}=\arccos{\left(\dfrac{L^2+(\|\mathbf{o}_{4t+3}-\mathbf{o}_{4t}\|/2)^2-(L/2)^2}{2L(\|\mathbf{o}_{4t+3}-\mathbf{o}_{4t}\|/2)}\right)}.
    \end{split}
\end{align}
Eqs.~\eqref{eq:turning-angle-1} and \eqref{eq:turning-angle-2} provide a reasonable setting to reproduce the pattern shown in Fig.~3(F).
Specifically, by inputting a straight polyline with segment length $L = 1$ and substituting Eqs.~\eqref{eq:turning-angle-1} and \eqref{eq:turning-angle-2} into Eq.~\eqref{eq:four-points} under the parameters $l = 1/3$ and $\phi^{*} = \phi_{0} = \pi/3$, we recover the pattern shown in Fig.~3(F) of the main text.
The central creases shown in Figs.~4(C) and 4(D) of the main text are obtained by applying the above procedure to an S-shaped input polyline with segment length $L = 1$, using the same parameters $l = 1/3$ and $\phi^{*} = \phi_{0} = \pi/3$.
The central crease of the Miura-ori strip shown in Fig.~4(E) is constructed from the same input polyline, with $l = 1/3$ and $\phi^{*} = \phi_{0} = 40^\circ$, by reversing the signs of the rotation angles $\psi_{t}$ and $\bar{\psi}_{t}$ in Eq.~\eqref{eq:four-points}.

\section*{Supplementary Movies}

\textbf{Movie S1.}
Folding motion of the strip shown in Fig.~3(A) of the main text. 
The evolution of the corresponding cobweb plot and phase-space trajectory at each time step is also shown.

\textbf{Movie S2.}
Folding motion of the strip shown in Fig.~3(B) of the main text. 
The evolution of the corresponding cobweb plot and phase-space trajectory at each time step is also shown.

\textbf{Movie S3.}
Folding motion of the strip shown in Fig.~3(C) of the main text. 
The evolution of the corresponding cobweb plot and phase-space trajectory at each time step is also shown.

\textbf{Movie S4.}
Folding motion of the strip shown in Fig.~3(D) of the main text. 
The evolution of the corresponding cobweb plot and phase-space trajectory at each time step is also shown.

\textbf{Movie S5.}
Folding motion of the strip shown in Fig.~3(E) of the main text. 
The evolution of the corresponding cobweb plot and phase-space trajectory at each time step is also shown.

\textbf{Movie S6.}
Folding motion of the strip shown in Fig.~3(F) of the main text. 
The evolution of the corresponding cobweb plot and phase-space trajectory at each time step is also shown.

\textbf{Movie S7.}
Folding motion of the strip shown in Fig.~4(D) of the main text.

\textbf{Movie S8.}
Folding motion of the strip shown in Fig.~4(E) of the main text.

\textbf{Movie S9.}
Folding motion of the thickness accommodated strip using the method shown in Fig.~5(A) of the main text.

\textbf{Movie S10.}
Domino-like sequential deployment of the 3D-printed physical prototype shown in Fig.~5(D) of the main text.